\documentclass[makeidx]{amsart}
\usepackage{amsmath}
\usepackage{amssymb}
\usepackage{epsfig}
\usepackage{hyperref}
\newtheorem{Proposition}{Proposition}[section]

  \newtheorem{Lemma}[Proposition]{Lemma}
  
  \newtheorem{Theorem}{Theorem}[section]
 
 \newtheorem{Definition}[Proposition]{Definition}
 \newtheorem{Note}[Proposition]{Note}

\newcommand {\z}{{\noindent}}
\def\blackslug{\hbox{\hskip 1pt \vrule width 4pt height 8pt depth 1.5pt
\hskip 1pt}}
\def\qed{\quad\blackslug\lower 8.5pt\null\par}
 
\def\CC{\mathbb{C}}
 \def\RR{\mathbb{R}}
 \def\NN{\mathbb{N}}

\def\Re{\mathrm{Re}}
\def\Im{\mathrm{Im}}

\makeindex
\begin{document}

\author{O. Costin$^1$ and X. Xia$^1$ }\title{Global
  reconstruction of analytic functions from local expansions}
\gdef\shortauthors{O.  Costin and X. Xia} \gdef\shorttitle{Global analytic
  reconstruction} \thanks{$1$.  Department of Mathematics, Ohio State
  University, 231 W. 18th Ave., Columbus, OH, 43210.}
  \maketitle

\bigskip

\begin{abstract}
  A new summation method is introduced to convert a relatively wide family of infinite sums and local expansions into integrals. The integral
  representations yield global information such as analytic
  continuability, position of singularities, asymptotics for large values of
  the variable and asymptotic location of zeros.  There is a duality between
  the global analytic structure of the reconstructed function and the
  properties of the coefficients as a function of their index.  Borel
  summability of a class of divergent series follow as a byproduct.

\end{abstract}
\tableofcontents

\setcounter{equation}{0}\section{Introduction}

Finding the global behavior of an analytic function in terms of its Taylor
coefficients is a notoriously difficult problem.  In fact, there cannot exist
a general solution to this problem, since many undecidable questions can be
quite readily reformulated in such terms.  

However, with some restrictions on the analytic functions, or correspondingly
on their Taylor coefficients $\{c_k\}_{k\in \NN}$ an effective and global
reconstruction becomes possible.  In the present paper we show that for the
class $\mathcal{M}$ of functions analytic in the complex plane with finitely
many cuts and with algebraic behavior at infinity (see Definition \ref{Def21}
below) global information is, perhaps surprisingly, contained in an effective
way in $\{c_k\}_{k\in \NN}$. In fact, the $c_k$'s  of such functions have
distinctive asymptotic features that allow for a global integral
representation of the associated function using a new summation method that
we introduce.  The position and type of singularities and the asymptotic
behavior of $f$ can be found in a practical way.  There is a duality between
the properties of the coefficients and the global structure, for instance
monodromy, of the reconstructed function.

More precisely, a function $f$ given by a series $f(z)=\sum_{k=0}^\infty c_k
z^k$ convergent for small $z$ belongs to $\mathcal{M}$ {\em if and only if}
the coefficients $c_k$ admit asymptotic representations in the form of
combinations of exponentials and Borel sums of (typically divergent) series.
In the language of generalized (EB, or \'Ecalle-Borel) summation the latter
condition is that the $c_k$'s have Borel summable transseries. Knowledge of
transseries is however not needed for converting a Taylor series into an
integral expression valid globally. However, the connection with transseries
is interesting and we describe it together with some key elements of
\'Ecalle-Borel summation in the Appendix \S\ref{EB}. This connection also shows that this class of Taylor
coefficients is fairly general. Further generalizations are discussed in \S\ref{Ge} and \S\ref{EB}.

The techniques we introduce also allow for a closed form
representation of functions for which the Taylor coefficients are known
explicitly such as those in \eqref{eq:ln2}.  We recently used this approach to analyze a class of linear PDEs with variable coefficients,
\cite{CHT}.

In particular, if $c_k=\varphi(k)$ where the function $\varphi$,  defined in the right
half plane,  is inverse Laplace transformable and its inverse Laplace
transform $\mathcal{L}^{-1} \varphi$ can be calculated in closed form, the function
$f$ has integral representations in terms of $\varphi$.   Some explicit examples we use for illustration are
\begin{equation}
  \label{eq:ln2}
 c_k^{[1]}=\frac{1}{(k+a)^b};  \ \ c_k^{[2]} =\frac{1}{k^b+\ln
    k},\ \  c_k^{[3]}= \frac{1}{k^{k+1}}, \
  c_k^{[4]}=e^{\sqrt{k}}, \ (a,b>0)
\end{equation}
We find that
\begin{equation}
  \label{eq:sum0}
  f_1(z)=:\sum_{k=1}^\infty c_k^{[1]}z^k= \frac{ z}{\Gamma(b)}\int_0^{\infty}\frac{[\ln(1+t)]^{b-1} dt}{(1+t)^{a+1}(t-(z-1))}
\end{equation}
On the first Riemann sheet $f_1$ has only one singularity, at $z=1$, of
logarithmic type, and $f_1=o(z)$ for large $z$.  General Riemann surface
information and monodromy follow straightforwardly from \eqref{eq:sum0}.  A
similar complex analytic structure is shared by $f_2=\sum_{k=1}^\infty
c_k^{[2]}z^k$, which has one singularity at $z=1$ where it is analytic in
$\ln(1-z)$ and $(1-z)$; more precisely,
\begin{equation}
  \label{intlog}
 f_2(z)=-\frac{1}{2\pi i} \frac{\ln \ln
   z}{z}\oint_{0}^{\infty}\frac{e^{-u\ln(z)}}{(-u)^{b}+\ln (-u)}du +E(z)
\end{equation}
see Definition \ref{def2}, where $E$ is entire.

  The function $f_3(z)=\sum_{k=1}^\infty
c_k^{[3]}z^k$ is entire; questions answered regard
say the behavior for large negative $z$ or the asymptotic location of
zeros. It will follow that $f_3$ can be written as
\begin{equation}
  \label{eq:invl}
 f_3(z)=\int_0^{\infty}(1+u)^{-1}G(\ln(1+u))\left[\exp\left({\frac{ze^{-1}}{1+u}}\right)-1\right]du
\end{equation}
where $G(p)=s'_2(1+p)-s'_1(1+p)$ and $s_{1,2}$ are two branches of the
functional inverse of $s-\ln s$, cf. \S\ref{entire}.  Using the integral
representation of $f_{3}$, its behavior for large $z$ can be obtained from
(\ref{eq:invl}) by standard asymptotics methods; in particular, for large
negative $z$, $f_3$ behaves like a constant plus $z^{-1/2}e^{-z/e}$ times a
factorially divergent series (whose terms can be calculated).

For 
$ c_k^{[4]}$ we find 
\begin{equation}
  \label{ff4}
f_4=\sum_{k=1}^\infty c_k^{[4]}z^k =-\frac{1}{4\sqrt{\pi}}\int_{C_1}p^{-3/2} \frac{e^{p+\frac{1}{4p}}}{e^p-z} e^{-np}dp
\end{equation}
where $C_1$ starts along $\RR^+$, loops clockwise once around the origin
and ends up at $+\infty$.

We also show that Borel summation of divergent series or transseries of
resurgent functions with finitely many Borel-plane singularities, as well as
the Abel-Plana version of the Euler-Maclaurin summation formula (see also
\cite{Stavros1}) can be derived by a natural extension of our analysis.

A separate category is represented by lacunary series. Their coefficients do
not satisfy our assumption; however a slightly different approach allows for a
detailed study of the associated functions as the natural boundary is
approached, \cite{Advances}.

\section{Main results}

\subsection{Global description from local expansions}$ $\\
$ $\ \ \ {\em A first class of problems} is finding the location and type of
singularities in $\CC$ and the behavior for large values of the variable of
functions given by series with finite radius of convergence (Theorem
\ref{T1}), such as those in (\ref{eq:ln2}).

{\em The second class of problems} amenable to the techniques presented
concerns the behavior at infinity (growth, decay, asymptotic location of
zeros etc.) of entire functions presented as Taylor series (Theorem \ref{T2}).

{\em The third question} is essentially the converse of the two above: given a
function that has analytic continuation on some Riemann surface, how is this
reflected on $c_k$? (Theorem \ref{T1}.)

{\em The third type of class of problems} is to determine Borel summability of
series with zero radius of convergence such as
\begin{equation}
  \label{eq:entire}
  \tilde{f}_4=\sum_{n=0}^{\infty}n^{n+1} z^n
\end{equation}
in which
the coefficients of the series are analyzable (Theorem \ref{T3}).

\begin{Definition}\label{Def21}
  {\rm Let $\{a_j : 1\leq j \leq N\}$ be a set of nonzero complex numbers with distinct arguments. Let $\mathcal{M}$ consist of the functions algebraically bounded
at $\infty$ and  analytic in $\CC\setminus \bigcup_{j=1}^N
 \{a_j t: t\geq 1\}$ .  By taking sufficiently many derivatives we can assume that
  $f\in\mathcal{M}'=\{f\in \mathcal{M}:f=o(z)\ \text{as}\
  z\to\infty\}$. }\\
  {\rm This is one of the simplest settings often occurring in applications. We can see later from the proof that the approach is more general.}
\end{Definition}
\begin{Definition}
\label{def2}
{\rm Assume $g(s)$ is analytic in $U_\delta \backslash [0,\infty) $ for some $\delta > 0$, where $U_\delta = \{z : |\Im (z)| \leq \delta, \Re(z) \geq -\delta \}$ and $g(s) \to 0$ uniformly in $U_\delta$, as $\Re(s) \to \infty$.
Assume $\epsilon \leq \delta$. We define $\Gamma_\epsilon$ to be the contour around $\RR^+$ consisting of two rays $l_{1,\epsilon}, l_{2,\epsilon}$ and a semicircle $\gamma_{\epsilon}$, where $l_{1,\epsilon} = \{x - \epsilon i: a\in[0,\infty)\}$ oriented towards the left, $l_{2,\epsilon} = \{x + \epsilon i: a\in[0,\infty)\}$ oriented towards the right; $\gamma_{\epsilon}$ is the left semicircle centered at origin oriented clockwise.
We denote 
\begin{equation}
  \label{eq:sg}
  \oint_0^{\infty}g(s)ds 
\end{equation}
the integral of $g$ over $\Gamma_\epsilon$. Since $g(s)$ vanishes at $\infty$, the integral is independent of the choice of $\epsilon$ as long as it is small enough.}
\end{Definition}
\begin{Note}
  {\rm A representation of the form \eqref{eq:sg} exists for Laplace
    transforms $h(t)= \int_0^{\infty}(s-t)^{-1}H(s)ds$ with $H$ analytic at
    zero, for instance  $h(t)=- (2\pi i)^{-1}
    \oint_0^{\infty}(s-t)^{-1}H(s)\ln s\,ds$.}
\end{Note}
\begin{Note}
  {\rm Consider the composition of $g$ with $s \to \ln (1+s)$, the branch cut of which is chosen to be $(-\infty, -1]$. If $g$ is analytic in $U_\delta \setminus [0,\infty)$, then $g(\ln(1+s))$ is analytic in the set $-1+\exp(U_\delta \setminus [0,\infty))$. If in addition we have the decay condition $g(\ln(1+s)) = o(s^{-(2-\alpha)})$ for some $\alpha \in (0,1)$ as $|s| \to \infty$, then there exists a $\tilde{\delta}$ small enough such that $U_{\tilde{\delta}} \subseteq  -1+\exp(U_\delta \setminus [0,\infty))$. It is easy to see from the decay condition that}
  \begin{equation*}
  \int_{\Gamma_{\tilde{\delta}}} g(\ln(1+s)) ds = \int_{\exp(\Gamma_\delta) -1} g(\ln(1+s)) ds
  \end{equation*}
  {\rm and so we can also choose the homotopic contour $\exp(\Gamma_\delta) -1$ instead of $\Gamma_{\tilde{\delta}}$ in the definition of $\oint_{0}^{\infty} g(\ln(1+s)) ds$. See Figure 1 for the contours.}
\end{Note}
\begin{figure}
\epsfig{file=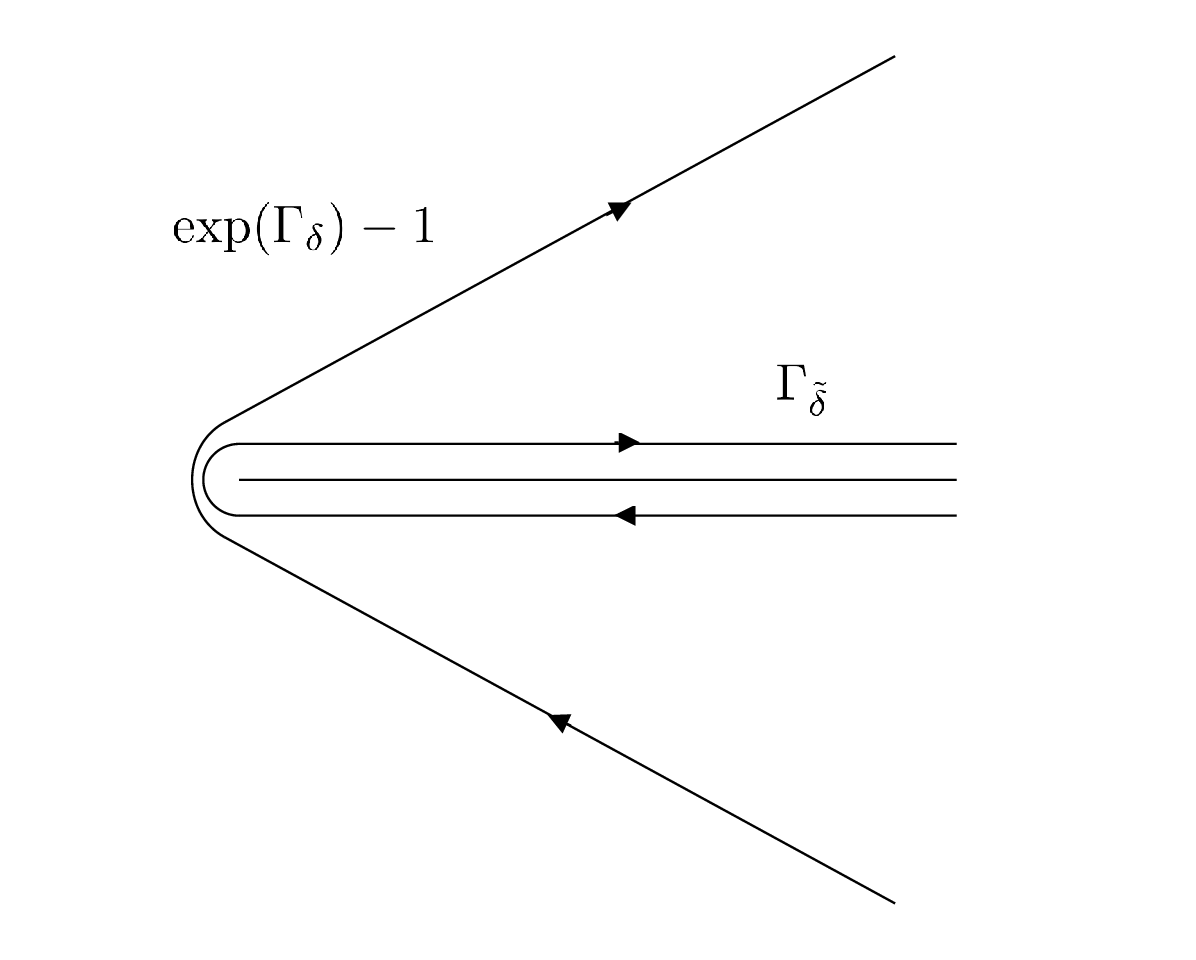, width=7cm,height=5cm}
\caption{}
\end{figure}

While providing integral formulae in terms of functions with known
singularities which are often rather explicit, the following result can also
be interpreted as a {\em duality of resurgence}.  \footnote{After developing these methods,
  it has been brought to our attention that a duality between resurgent
  functions and resurgent Taylor coefficients has been noted in an unpublished
  manuscript by \'Ecalle.}.
\begin{Theorem}\label{T1}
  (i) Assume that $f(z)=\sum_{k=0}^{\infty}c_k z^k$ is a series with positive, finite radius of convergence,  with $c_k$ having Borel sum-like representations of the form
\begin{equation}
  \label{eq:trs1}
  c_k=\sum_{j=1}^N a_j^{-k}\oint_0^{\infty}e^{-kp}F_j(p)dp \hspace{0.4in} (k \ge 1)
\end{equation}
(\ref{eq:trs1}) with $a_j$ as in Definition 2.1, $F_j$ analytic in $U_\delta \backslash [0,\infty) $ for some $\delta > 0$ and algebraically bounded at $\infty$. Then, $f$ is given by
\begin{equation}
  \label{eq:recon3}
 f(z)=f(0)+z\oint_0^{\infty}\sum_{j=1}^N \frac{F_j(\ln (1+s))ds}{(1+s)((1+s) a_j - z)}
\end{equation}

\z (ii) Furthermore, $f \in \mathcal{M}'$.The behavior of $f$ at $a_j$ and is of the same type as the behavior of $F_j(\ln(1+s))$ at $0$. More precisely, for small $z \notin [0,\infty)$,
\begin{equation}
\label{eq:singtype}
f(a_j (z+1)) = 2\pi i F_j (\ln (1+z)) + G(z)
\end{equation}
where $G(z)$ is analytic at 0.

\z (iii) Conversely, assume $f \in \mathcal{M}'$, and has finitely many singularities located at $\{a_j t_{j,l} \}$, $(1 \leq j \leq N, 1\leq l)$, with $1 =t _{j,1}$ and $ t _{j,l} < t_{j,l+1}$ for all $j, l$. Let
$c_k=f^{(k)}(0)/k!$; then $c_k$ have Borel sum-like representations of the form
\begin{equation}
  \label{eq:trs12}
  c_k=\frac{1}{2\pi i}\sum_{j=1}^N (a_j )^{-k}\oint_0^{\infty}e^{-ks}\ f( a_{j} e^s)ds,\  k \geq 1
\end{equation}

\end{Theorem}

 The behavior at $a_j$ and
at $\infty$ will follow from the proof.

\bigskip

\begin{figure}
\epsfig{file=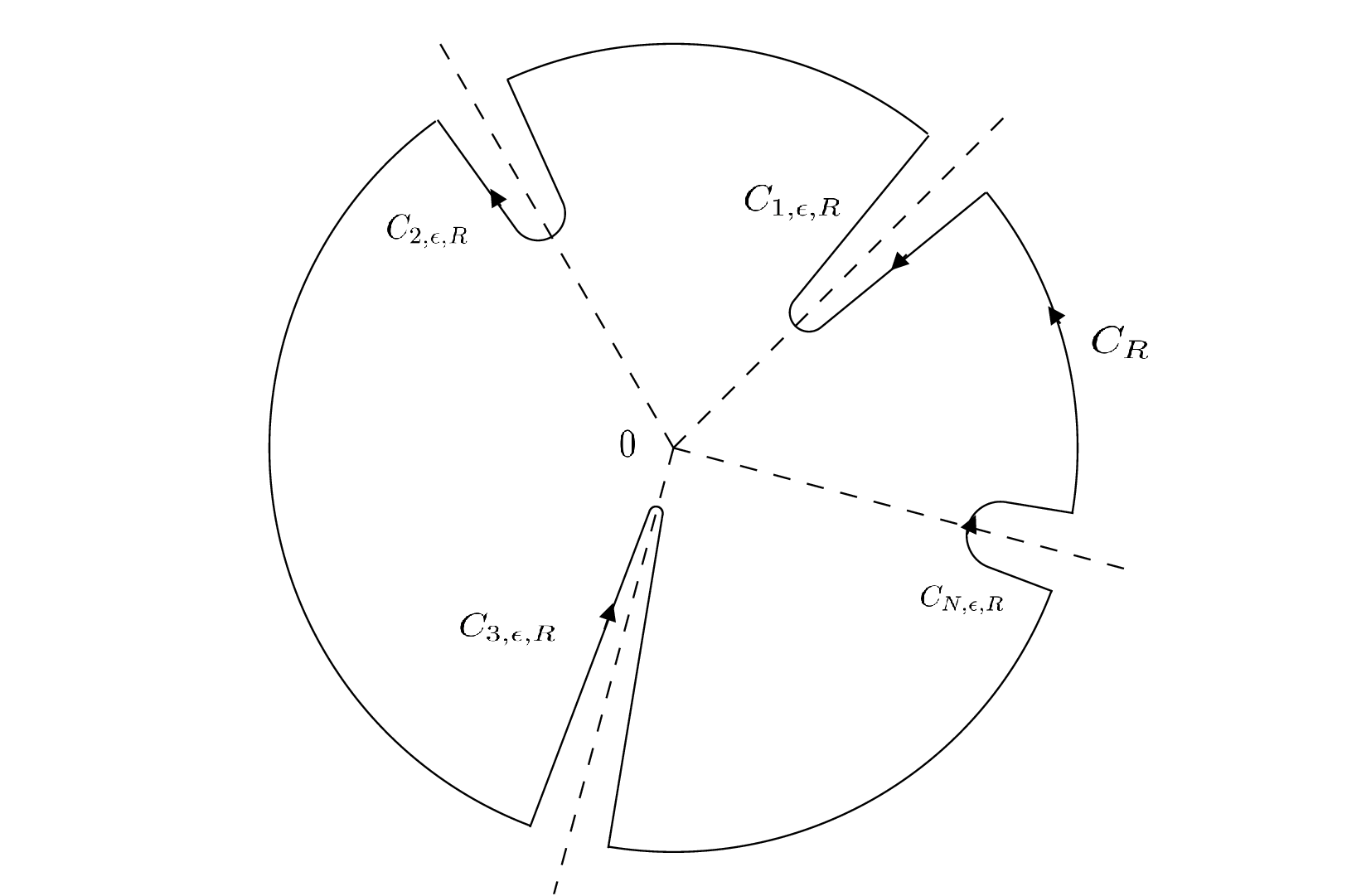, width=12cm,height=8cm}
\caption{Singularities of $f$, cuts, direction of integration and Cauchy
  contour deformation.}
\end{figure}

As it will be clear from the proofs, the method and results would apply, with
minor adaptations to functions of several complex variables.

\subsection{Entire functions}\label{entire}

 We
restrict the analysis to entire functions of exponential order one, with
complete information on the Taylor coefficients. Such functions include of
course the exponential itself, or expressions such as $f_3$. It is useful to
start with $f_3$ as an example. The analysis is brought to the case in
\S\ref{CaseI} by first taking a Laplace transform. Note that
\begin{equation}
  \label{eq:eq41}
  \int_0^{\infty}e^{-xz}f(z)dz=\frac{1}{x}\sum_{n=1}^\infty\frac{n!}{n^{n+1}x^{n}}
\end{equation}
The study of entire functions of exponential order one likely involves the
factorial, and then a Borel summed representation of the Stirling formula is
needed.
\begin{Theorem}\label{T2}
  Assume that the entire function $f$ is given by
\begin{equation}
  \label{eq:enti22}
 f(z)= \sum_{k=1}^\infty \frac{c_k z^k }{k!}
\end{equation}
with $c_k$ as in Theorem~\ref{T1}. Then,
\begin{equation}
  \label{eq:recon31}
 f(z)=
  \oint_0^{\infty}\sum_{j=1}^N \left[\left(e^{\frac{z}{s+a_j}}-1\right)
\frac{F_j(\ln (1+s))}{(1+s)}\right]ds
\end{equation}
\end{Theorem}
As in the simple example, the behavior at infinity follows from the integral
representation by classical means.
\subsection{Borel summation}
We obtain from Theorem~\ref{T1}, in the same way as above, the
following.
\begin{Theorem}\label{T3}
  Consider the formal power series
\begin{equation}
  \label{eq:enti23}
 \tilde{f}(z)= \sum_{k=1}^\infty {c_k k! z^{k+1}}
\end{equation}
with coefficients $c_k$ as in Theorem~\ref{T1}.
Then the series \eqref{eq:enti23} is (generalized)  Borel summable to
\begin{multline}
  \label{eq:bs1}
\int_0^{\infty}dp e^{-p/z}p\sum_{j=1}^N\oint_0^{\infty}\frac{F_j(\ln(1+s))}{(1+s)(a_js+a_j-p)}ds\\=
   -\sum_{j=1}^N\oint_0^{\infty}\frac{F_j(\ln (1+s))}{1+s}\left(z-a_j(s+1)e^{-\frac{a_j(s+1)}{z}}\mathrm{Ei}\left(\frac{a_j(s+1)}{z}\right)\right)ds
\end{multline}
\end{Theorem}
The proof proceeds as in the previous sections, taking now a Borel
 followed by Laplace transform.
\subsection{Other applications; the examples in the introduction}\label{Examples}
\subsubsection{Other growth rates}\label{Growthrates}

Other growth rates can be accomodated, for instance by analytic continuation. 
We have for positive $\gamma$,
\begin{equation}
  \label{eq:eqinvl}
  e^{-\gamma\sqrt{
      n}}=\frac{\gamma}{2\sqrt{\pi}}\int_0^{\infty}p^{-3/2}e^{-\frac{\gamma^2}{4p}}e^{-np}dp
\end{equation}
which can be analytically continued in $\gamma$. We note first that the
contour cannot be, for this function, detached from zero. Instead, we keep one
endpoint at infinity and, near the
origin, simultaneously rotate $\gamma$ and $p$  to maintain $-\gamma/p$ real
and negative. We get 
\begin{equation}
  \label{eq:eqerfc}
  e^{\sqrt{n}}=-\frac{1}{4\sqrt{\pi}}\int_{C_1} p^{-3/2} e^{\frac{1}{4p}}e^{-np}dp
\end{equation}
 and \eqref{ff4} follows.  In particular, 
\begin{equation}
  \label{eq:limit1}
  \lim_{z\to -1^+}\sum_{n=1}^{\infty}e^{\sqrt{n}}z^n=
  -\frac{1}{4\sqrt{\pi}}\int_{C_1}\frac{e^{1/p}dp}{p^{-3/2}(e^p+1)}
\end{equation}
 The sum \eqref{eq:limit1} is  unwieldy numerically,
while the integral \eqref{eq:eqinvl} can be evaluated accurately by standard
means. In a similar way we get
\begin{equation}
  \label{eq:eqsum}
  \sum_{k=0}^{\infty}\frac{e^{i\sqrt{k}}}{k^{a}}=-\frac{\gamma
    2^{a-1/2}}{\sqrt{\pi}}
\int_Cdp \frac{e^{-\frac{1}{8p}}U(2a+1/2;\frac{1}{\sqrt{2p}})}{p^{a-1}(e^p+1)}
\end{equation}
for $a>1/2$ for which the series converges.
Here $C$ is a contour starting along $\RR^-$, circling the origin clockwise
and ending up at $+\infty$ and $U$ is the parabolic cylinder function
\cite{Abramowitz}.  These sums are obtained in \S\ref{Growthrates}.

{\em The coefficients $c_k^{[1]}$ in \eqref{eq:ln2}}.
We have
\begin{equation}
  \label{eq:invlp3}
  \mathcal{L}^{-1}\left[\frac{1}{(n+a)^b}\right]=\Gamma(b)^{-1}p^{b-1}e^{-ap}
\end{equation}
The rest follows in the same way \eqref{ff4} was obtained, after changing
variables to $1+t=e^p$.

{\em The coefficients $c_k^{[2]}$}.
 We let $x=n$ and take the inverse
Laplace transform in $x$:
\begin{equation}
  \label{eq:dual1}
  G(p) = \frac{1}{2\pi i}\int_{1-i\infty}^{1+i\infty}\frac{e^{xp}}{x^{b}+\ln x}dx
\end{equation}
where the contour can be bent backwards for $p \in \RR^+$, to hang around $\RR^-$. Then, with
the change of variable $x=-u$  (\ref{eq:dual1}) becomes
\begin{equation}
  \label{eq:dual2}
 G(p)= \frac{1}{2 \pi i}\oint_{0}^{\infty}\frac{e^{-up}}{(-u)^{b}+\ln (-u)}du
\end{equation}
and thus
\begin{equation*}
c_k = (\mathcal{L} G) (k) = \int_0^{\infty} G (p) e^{-kp}{dp }  = \oint_0^{\infty} \left[ \frac{-G(p) \ln p}{ 2 \pi i} \right] e^{-kp}{dp}
\end{equation*}
We see that $F_1(p) = (- G(p) \ln p)/{2 \pi i}$ and by Theorem 2.1
\begin{equation*}
f_2 (z) = z \oint_0^{\infty} \frac {\tilde{G}(s)}{s - (z-1)} {ds}
\end{equation*}
where $\tilde {G}(s) = F_1 (\ln(1+s))/(1+s)$.
Hence the singularity of $f_2(z)$, at one, according to (\ref{eq:sing0}) is
that of $ \phi(z) =  2 \pi i \tilde{G}(z - 1)$, as in
(\ref{intlog}).

The example of {\em the coefficients $c_k^{[3]}$} is an application of Theorem
\ref{T2}, after calculations similar to the ones above.

{\em The coefficients $c_k^{[4]}$} were treated at the beginning of this
  section.

{\em Another example} is provided by the log of the Gamma function,
$\ln\Gamma(n)=\sum_{k=1}^n \ln k$. It is convenient to first subtract out the leading
behavior of the sum to arrange that the summand is inverse Laplace
transformable. With 
$g_n=\ln\Gamma(n)-(n\ln n-n-\frac{1}{2}\ln n)$  we get 
\begin{multline}
  g_n=\sum_{1}^n
  \Big[1-\Big(\frac{1}{2}+n\Big)\ln\Big(1+\frac{1}{n}\Big)\Big]=\sum_{1}^n\int_0^\infty
  e^{-np} \frac{ 1-\frac{p}{2}-(\frac{p}{2}+1)e^{-p}}{p^2}dp
\end{multline}
where $\mathcal{L}^{-1}$ of the summand in the middle term is most easily
obtained by noting that its second derivative is a rational function. 
Summing as usual $e^{-np}$ we get
\begin{equation}
  \label{eq:lngamma}
  \ln\Gamma(n)=n(\ln n -1)-\frac{1}{2}\ln n +\frac{1}{2}\ln(2\pi)+\int_0^{\infty}\frac{\displaystyle
  1-\frac{p}{2}-\Big(\frac{p}{2}+1\Big)e^{-p}}{p^2(e^{-p}-1)}e^{-np}dp
\end{equation}

Obviously, if the behavior of the coefficients is of the form $A^kc_k$ where
$c_k$ satisfies
the conditions in the paper, one simply changes the independent variable to  $z'=Az$.

\subsection{The Gamma function and Borel summed Stirling formula}\label{S115}
We have
\begin{multline}
  n!=\int_0^\infty t^ne^{-t}dt=n^{n+1}\int_0^\infty e^{-n(s-\ln s)}ds\\=
  n^{n+1}\int_0^1 e^{-n(s-\ln s)}ds+ n^{n+1}\int_1^\infty e^{-n(s-\ln s)}ds
\end{multline}

\z On $(0,1)$ and $(1,\infty)$ separately, the function $s-\ln(s)$ is
monotonic and we may write, after inverting $s-\ln(s)=t$ on the two
intervals to get $s_{1,2}=s_{1,2}(t)$ \footnote{The functions $s_{1,2}$ are
  given by branches of $-W(-e^t)$, where $W$ is the Lambert function.},
\begin{equation}\label{e25}
 n!= n^{n+1}\int_1^{\infty}e^{-nt}(s'_2(t)-s'_1(t))dt=n^{n+1}e^{-n}
\int_0^{\infty}e^{-np}G(p)dp
\end{equation}
\z where $G(p)=s_2'(1+p)-s_1'(1+p)$. From the definition it follows that $G$
is bounded at infinity and $p^{1/2}G$ is analytic in $p$ at $p=0$.
Using now (\ref{e25}) and Theorem~\ref{T1} in (\ref{eq:eq41}) we get
\begin{equation}
  \label{eq:lapl41}
  \int_0^{\infty}e^{-xz}f(z)dz=\frac{1}{ x^2}\int_0^{\infty}\frac{G(\ln(1+t))}{(te+(e-x^{-1}))(t+1)}dt
\end{equation}
Upon taking the inverse Laplace transform we obtain (\ref{eq:invl}).

\section{Proof of Theorem \ref{T1}}
\z If $f\in\mathcal{M}'$ we write the Taylor coefficients in the form
\begin{equation}
      \label{eq:cf1}
    c_k=\frac{1}{2\pi i}\oint \frac{dp f(p)}{p^{k+1}} \hspace{0.4in} (k\geq 1)
  \end{equation}
  where the contour of integration is a small circle of radius $r$ around the
  origin. We attempt to increase $r$ without bound. In the process, the contour
  will hang around the singularities of $f$ as shown in Figure 2. The
  integrals converge by the decay assumptions and the contribution of the arcs
  at large $r$ vanish, since $f(z) = o(z)$ as $z \to \infty$.

  To be more precise, let $C_{j,\epsilon}$ be the part of  the image of $\Gamma_{\epsilon}$ under the mapping $s \to a_j e^s $, let $C_{j,\epsilon, R}$ be the part of $C_{j,\epsilon}$ inside the disk $|s| \leq R$, and $C_R$ be the part of the contour on $|s| = R$. Then for $\epsilon$ small enough and $R$ large enough we have
\begin{equation*}
c_k=\frac{1}{2\pi i}\oint \frac{dp f(p)}{p^{k+1}} = \frac{1}{2\pi i} \left( \sum_{j = 1}^{N} \int_{C_{j, \epsilon, R}} \frac{dp f(p)}{p^{k+1} }+ \int_{C_R} \frac{dp f(p)}{p^{k+1} }\right)
\end{equation*}
By the change of variable $p = a_j e^s$ and letting $R \to \infty$ we get
\begin{equation*}
\int_{a_j C_{\epsilon, R}} \frac{dp f(p)}{p^{k+1} } = \int_{\Gamma_{\epsilon}} \frac{a_j e^s ds f(a_j e^s)}{(a_j e^s)^{k+1} } = \oint_{0}^{\infty} {a_j}^{-k} e^{-ks} f(a_j e^s) ds
\end{equation*}  
and the integral over $C_R$ vanishes as $R \to \infty$ by decay condition for $k \geq 1$.

\bigskip

  \z In the opposite direction, first let $\epsilon$ small enough so that for all $j$ and $k = 1$
  \begin{equation}
  \label{eq:path}
  \oint_{0}^{\infty}{e^{-kp}F_j(p) dp} = \int_{\Gamma_{\epsilon}}{e^{-kp}F_j(p) dp} 
  \end{equation}
 Then for all $1\leq j \leq N$, $k \geq 1$ (\ref{eq:path}) is true. Also let $z$ be small so that 
  \begin{equation}
  \label{eq:est1}
  |{a_j}^{-1}e^{-p} z| \leq \delta^2 < 1
  \end{equation}
  for all $j$ and $p \in \Gamma_\epsilon$.

Then, by the dominated convergence theorem (which applies in this case, see (\ref{eq:estimate})) we have
\begin{align}
  \label{eq:converse}
  &f(z) - f(0) 
  = \sum_{k=1}^{\infty}c_k z^k
  =\sum_{k=1}^{\infty}\left({\sum_{j=1}^N a_j^{-k}\int_{\Gamma_\epsilon}e^{-kp}F_j(p)dp} \right) z^k\\ \notag
&=\sum_{k=1}^{\infty}\int_{\Gamma_\epsilon} {\left( \sum_{j=1}^N (a_j^{-1}e^{-p} z)^kF_j(p) \right)} dp = \int_{\Gamma_\epsilon} \sum_{k=1}^{\infty} {\left( \sum_{j=1}^N (a_j^{-1}e^{-p} z)^kF_j(p) \right)} dp\\ \notag
&= \int_{\Gamma_\epsilon} \sum_{j=1}^{N} \left( \sum_{k=1}^{\infty}  (a_j^{-1}e^{-p} z)^k \right) F_j(p) dp = \int_{\Gamma_\epsilon} \sum_{j=1}^{N} \left( \frac{a_j^{-1}e^{-p} z}{1 - a_j^{-1}e^{-p} z} \right) F_j(p) dp\\ \notag
&=\sum_{j=1}^{N} \int_{\Gamma_\epsilon} \frac{a_j^{-1}e^{-p} z}{1 - a_j^{-1}e^{-p} z} F_j(p) dp = \sum_{j=1}^N z\int_{\exp(\Gamma_\epsilon) -1}\frac{F_j(\ln (1+s))ds}{(1+s)(s a_j+a_j-z)}\\ \notag
&=\sum_{j=1}^N 
  z\oint_0^{\infty}\frac{F_j(\ln (1+s))ds}{(1+s)(s a_j+a_j-z)}
\end{align}
as stated.  The fourth equality holds because we have, in view of \eqref{eq:est1},
\begin{multline}
\label{eq:estimate}
\left| \sum_{j=1}^N (a_j^{-1}e^{-p} z)^kF_j(p)\right| \leq \sum_{j=1}^N |a_j^{-1} e^{-p} z|^{k/2} \left( |a_j^{-1} e^{-p} z|^{k/2} |F_j(p)|\right) \\
\leq \sum_{j=1}^N \delta^k \left( |a_j^{-1} e^{-p} z|^{k/2} |F_j(p)|\right)
\end{multline}
For each $j$, $|a_j^{-1} e^{-p} z|^{k/2} F_j(p)$ is integrable over $\Gamma_\epsilon$ since $F_j$ is algebraically bounded at $\infty$. Hence we may interchange the order of integration and summation over $k$.\\
Furthermore, $f(z)$ can be analytically continued to a function in $\mathcal{M}'$, see  \eqref{eq:recon3}. Let $C$ be an arbitrary closed curve in $\CC\setminus \bigcup_{j=1}^N \{a_j t: t\geq 1\} $ and choose $\epsilon$ small enough so that for each $j$, ${\rm dist}(C, (\Gamma_\epsilon+1)a_j) > 0$. Since $F_j$ is algebraically bounded, we have
\begin{equation*}
\oint_0^{\infty}\frac{F_j(\ln (1+s))ds}{(1+s)(s a_j+a_j-z)} = \int_{\Gamma_\epsilon}\frac{F_j(\ln (1+s))ds}{(1+s)(s a_j+a_j-z)}
\end{equation*}
now, $|(s+1)a_j - z|$ is bounded below for $s\in \Gamma_\epsilon$, and $F_j(\ln (1+s))/\left((1+s)(s a_j + a_j - z) \right)= o(s^{-(2-\alpha)})\, (s \to \infty, \alpha \in (0,1))$. Thus, we have, by Fubini,
\begin{multline*}
\int_{C} dz \oint_0^{\infty}\frac{F_j(\ln (1+s))ds}{(1+s)(s a_j+a_j-z)} = \int_C dz \int_{\Gamma_\epsilon}\frac{F_j(\ln (1+s))ds}{(1+s)(s a_j+a_j-z)} =  \\
\int_{\Gamma_\epsilon} ds \int_C\frac{F_j(\ln (1+s))dz}{(1+s)(s a_j+a_j-z)} = \int_{\Gamma_\epsilon}\frac{F_j(\ln (1+s))ds}{(1+s)}\int_C\frac{dz}{s a_j+a_j-z} = 0
\end{multline*}
Morera's Theorem implies that the integral in (\ref{eq:recon3}) is analytic in $\CC\setminus \bigcup_{j=1}^N \{a_j t: t\geq 1\} $. In other words, $f(z)$ can be analytically continued in $z$. \\
\begin{figure}
\epsfig{file=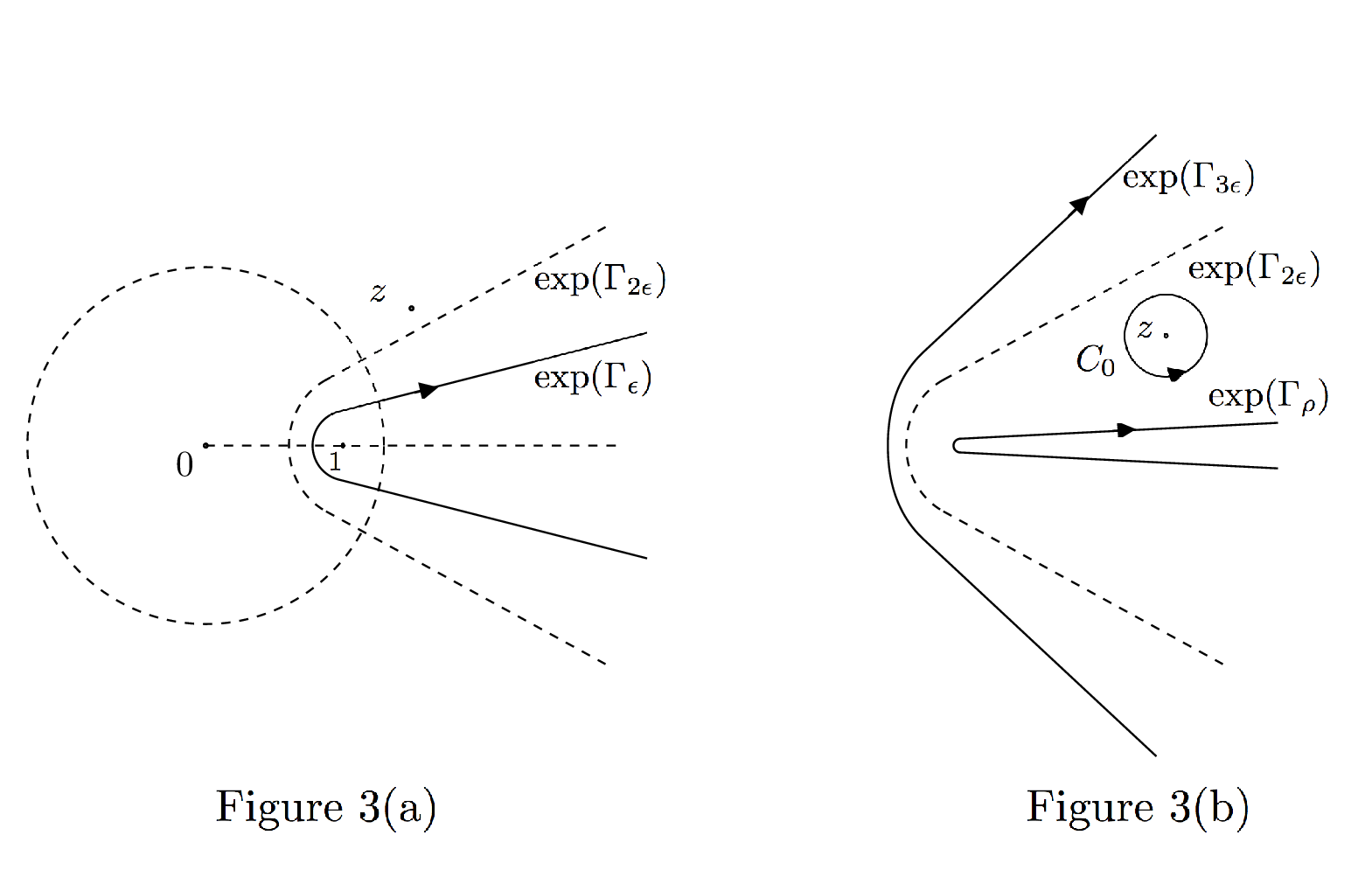, width=13cm,height=8cm}
\caption{Contours used in the proof of (\ref{eq:iii}). In Figure 3(a) $\arg|z| \in [2\epsilon, \pi]$ while in  3(b) $\arg|z| \in [0, 2\epsilon]$}
\end{figure}

Finally we show that $f(z) = o(z)$ and thus is in $\mathcal{M}'$. It suffices to show that  in the case $a_j = 1$,
\begin{equation}
\label{eq:iii}
\oint_{0}^{\infty} \frac{F_j(\ln(1+s)) ds}{(1+s)(1+s -z)} \to 0 \hspace{0.4 in}  (|z| \to \infty)
\end{equation}
Assume that $F_j$ is analytic in $U_{3\epsilon} \setminus [0,\infty)$. Note that since the integrand is $o(s^{-(2-\alpha)})$ at $\infty$, if $|\arg(z)| \in [2\epsilon, \pi]$ we may choose the contour to be $\exp(\Gamma_\epsilon) - 1$ as shown in Figure 3(a). Without loss of generality, assume $\arg(z) \in [2\epsilon, \pi]$, let $\beta = \arg(z) - \epsilon$, $\beta \in [\epsilon, \pi - \epsilon]$, and assume $|z| \geq 1$. Then we have
\begin{equation*}
\oint_{0}^{\infty} \frac{F_j(\ln(1+s)) ds}{(1+s)(1+s -z)} = \int_{\exp(\Gamma_\epsilon) - 1}\frac{F_j(\ln(1+s)) ds}{(1+s)(1+s -z)} =\int_{\exp(\Gamma_\epsilon)}\frac{F_j(\ln(t)) dt}{t(t -z)}
\end{equation*}
Notice that for $|t| \geq 1$, from the geometric interpretation or by direct calculation, we get
\begin{equation*}
|t-z|^2 \geq |t|^2 \sin^2 \beta
\geq |t|^2 \sin^2 \epsilon
\end{equation*}
Since $F_j$ is algebraically bounded at $\infty$, we have $F_j(\ln(t)) = o(t^{-(1-\alpha)})$ for any $\alpha \in (0,1)$, so by dominated convergence, as $|z| \to \infty$, we have
\begin{equation}
\label{eq:iiia}
\int_{\exp(\Gamma_\epsilon)}\frac{F_j(\ln(t)) dt}{t(t -z)} \to 0
\end{equation}
Let $\rho$ in $(0, \arg(z))$ recall that $|\arg(z)| \in (0,2\epsilon]$). Let $C_0$ be a positively oriented small circle around $z$ lying between the contours $\exp(\Gamma_\rho)$ and $\exp(\Gamma_{3\epsilon})$. By deforming  the contour $\exp(\Gamma_\rho)$ (see Figure 3(b)) we have
\begin{equation}
\label{eq:deform1}
\begin{split}
\oint_{0}^{\infty} \frac{F_j(\ln(1+s)) ds}{(1+s)(1+s -z)} = \int_{\exp(\Gamma_\rho)}\frac{F_j(\ln(t)) dt}{t(t -z)} = \left( \int_{\exp(\Gamma_{3\epsilon})} + \int_{C_0}\right)\frac{F_j(\ln(t)) dt}{t(t -z)} \\
=2\pi i \frac{F_j(\ln(z))}{z} + \int_{\exp(\Gamma_{3\epsilon})}\frac{F_j(\ln(t)) dt}{t(t -z)} 
\end{split}
\end{equation}
An analysis similar to the one leading to (\ref{eq:iiia}) shows that the integral over $\exp(\Gamma_{3\epsilon})$ vanishes as $|z| \to \infty$. The term before the last integral in \ref{eq:deform1} is simply $o(z^{-(1-\alpha)})$ for any $\alpha \in (0,1)$.
The nature of the singularities of $f$ is derived from \S\ref{CaseI}.
From (\ref{eq:converse}) we have 
\begin{equation*}
f(z) = f(0) + z \oint_{0}^{\infty} \sum_{j=1}^{N} \frac{F_j(\ln(1+s)) ds}{(1+s)((1+s)a_j -z)}
\end{equation*}
For $i \neq j$, if $z$ is close enough to $a_j$, and all the contours are chosen close enough to $[0,\infty)$, then
\begin{equation*}
\oint_{0}^{\infty} \frac{F_i(\ln(1+s)) ds}{(1+s)((1+s)a_j -z)} = G_i(z)
\end{equation*}
where $G_i(z)$ is analytic. Thus, for $z$ we have
\begin{align*}
&f((z+1)a_j)\\ &=f(0) + (z+1)a_j \oint_{0}^{\infty} \frac{F_j(\ln(1+s)) ds}{(1+s)((1+s)a_j -(z+1)a_j)}  + \sum_{i \neq j} G_i((z+1)a_j)\\
&= f(0) + (z+1)\oint_0^{\infty} \frac{F_j \ln(1+s) ds}{(1+s)(s-z)} +\sum_{i \neq j} G_i((z+1)a_j)\\
&= f(0) + (z+1)\left[ 2 \pi i \frac{F_j(\ln (1+z))}{1+z} + G(z)\right] +\sum_{i \neq j} G_i((z+1)a_j)\\
&= 2 \pi i F_j(\ln(1+z)) + \tilde{G}(z) 
\end{align*}
where $\tilde{G}(z)$ is analytic for small $z$, and hence (\ref{eq:singtype}) follows.
\subsubsection{Singularity formula. Duality.}\label{CaseI}
The problem of the type of singularities of the resummed series
reduces to finding the singularity type of a Hilbert-transform-like integral
of the form \eqref{eq:sg}.

The singularity of $g$ at $t=0$ is the same as the singularity of $G$ at $s=0$
as follows from a simple calculation:
\begin{Lemma}[Analytic structure at $t=0$.] Assume that  for some $\delta>0$, $G$ is analytic in $U_\delta \setminus [0,\infty)$ as in Definition \ref{def2} and $G(s) = o(s)$ at $\infty$. Then, 
\begin{equation}
    \label{eq:Plm1}
  g(t):= \oint_0^{\infty}\frac{G(s)}{s-t}ds=
2\pi i G(t)+G_2(t)
  \end{equation}
where $G_2(t)$ is analytic for small $t$.
\end{Lemma}
\begin{proof}
  The approach is similar to that of the proof of $f(z) = o(z)$ in the case where $\arg(z) \in [0,2\epsilon]$ in Theorem 2.1. We first note that for $t\in\CC\setminus [0,\infty)$ $g$ is analytic in
  $t$. To find the behavior of $g$ at $t=0$ we take a $t$ with $|t|=\epsilon$
  small enough outside the contour of integration. We next deform the contour
  around zero into a circle of radius $2\epsilon$ in the process collecting a
  residue
\begin{equation}
  \label{eq:sing0}
  2\pi i G(t)
\end{equation}
The new integral is manifestly analytic for $|t|<\epsilon$.

\end{proof}
\section{Generalizations}\label{Ge}

\begin{Note}
\begin{enumerate}

{\rm \item One can allow for
infinitely many singularities, under suitable estimates of their strengths.

\item Exponential growth of $F_j$  can be accomodated, 
provided a sufficient number of initial terms of the series in Theorem 2.1(i) are summed separately.

\item Several complex variables can be treated very similarly, see \cite{CHT}.
\item{Other types of decay/growth of coeefficients can be accomodated, cf.
\S\ref{Growthrates}}

}
\end{enumerate}
\end{Note}

\section{Appendix: Overview of transseries,  analyzable functions and EB summation}
\label{EB}
In the early 1980's, \'Ecalle discovered and extensively studied a broad class
of functions,  analyzable functions, closed under algebraic operations,
composition, function inversion, differentiation, integration and solution of
suitably restricted differential equations
\cite{Ecalle-book,Ecalle,Ecalle2,Ecalle3}.  They are described as
``transseries", generalized expansions obtained by {\sl closing power series under the same
operations}. Transseries are surprisingly easy to describe; roughly,
they are ordinal length, asymptotic expansions involving powers, iterated
exponentials and logs, with at most power-of-factorially growing coefficients.

In view of the closure of analyzable functions with respect to a wide class of operations,
reconstructing functions from series with arbitrary {\em analyzable}
coefficients would make the reconstruction likely applicable to series
occurring in problems involving any combination of these many operations.

The class of coefficients having \'Ecalle Borel (EB)-summable transseries is fairly
general. In particular  it is known \cite{Ecalle, Book} that solutions of linear or
nonlinear recurrence relations of arbitrary order with analytic coefficients
are EB-summable, \cite{Braaksma}.  Recurrence relations exist for instance
when the coefficients are obtained by solving differential equations by power
series.

This paper deals with analyzable coefficients having finitely many
singularities after a suitable Borel transform. The methods however
are open to substantial extension. In particular, we allow for general
singularities, while analyzable and resurgent functions have singularities
of a controlled type \cite{Ecalle}.

\subsection{Classical and generalized Borel summation} \label{Bsum} A series
$\tilde{f}=\sum_{k=1}^{\infty}c_k x^{-k}$ is Borel summable if its Borel
transform, {\em i.e.} the formal
inverse Laplace transform \footnote{$\displaystyle \sum_{k=1}^{\infty}c_k \frac{1}{2\pi i}\int_{1-i\infty}^{1+i\infty} e^{px}
  x^{-k}dx=\sum_{k=0}^{\infty}\frac{c_k p^{k-1}}{(k-1)!}=F(p)$}
  converges to a function $F$ analytic in a neighborhood of $\RR^+$, and $F$
  grows at most exponentially at infinity.  The Laplace transform of $F$ is by
  definition the Borel sum of $\tilde{f}$.  Since Borel summation is formally
  the identity, it is an extended isomorphism between functions and series,
  much as convergent Taylor series associate to their sums.

  However expansions occurring in applications are often not {\em classically}
  Borel summable, sometimes for the relatively manageable reason that the
  expansions are not simple integer power series, or often, more seriously,
  because $F$ is singular on $\RR^+$, as is the case of $\sum n!x^{-n-1}$
  where $F=(1-p)^{-1}$, or
  because $F$ grows superexponentially.

To address the latter difficulties, \'Ecalle defined averaging and cohesive
continuation to replace analytic continuation, and acceleration to deal with
superexponential growth \cite{Ecalle-book,Ecalle,Ecalle2,Ecalle3}.

We call \'Ecalle's technique \'Ecalle-Borel (EB) summability and ``EB
transform'' the inverse of EB summation. While it is an open, imprecisely
formulated, and in fact conceptually challenging question, whether EB summable
series are closed under all operations needed in analysis, general results
have been proved for ODEs, difference equations, PDEs, KAM resonant expansions
and other classes of problems \cite{Balser, Duke, Braaksma2,Poincare,OCSG}. EB
summability seems for now quite general.

A function is analyzable if it is an EB transform of a transseries. This
transseries is then unique \cite{Ecalle}. Then the EB transform is the mapping
that associates this unique transseries to the function. Simple examples of such transseries are
\begin{equation}
  \label{eq:dual}
  c_k=\frac{1}{\sqrt{k}};\ \ 
c_k=\frac{1}{k^b}\sum_{j=0}^{\infty} (-1)^k\frac{(\ln k)^j}{k^{j\pi}};\ \
c_k=\sum_{j=1}^{\infty}\frac{(-z)^j j!}{k^j}+2^{-k-1}\sum_{j=0}^{\infty}\frac{(-z)^j\Gamma(j-1/2)}{\sqrt{\pi}}
\end{equation}
The first two are convergent and correspond to the first two series in
(\ref{eq:ln2}); the last one is divergent but Borel summable.  

EB summation consists, in the simplest cases, in replacing the series in the
transseries by their Borel sum.

\subsubsection{ Borel summed version of  $1/n!$} 
We can use the following representation \cite{Abramowitz}
\begin{equation}
  \label{eq:eqg}
  \frac{1}{\Gamma(z)}=-\frac{ie^{-\pi i
      z}}{2\pi}\oint_0^{\infty}s^{-z}e^{-s}ds=-\frac{ie^{-\pi i
      z}z^{-z}}{2\pi}\oint_0^{\infty}s^{-z}e^{-zs}ds
\end{equation}
with our convention of contour integration. From here, one can proceed 
as in \S\ref{S115}.

\subsection{Acknowledgments}  M. Kruskal introduced OC to questions of the
type addressed in the paper; OC is grateful to S.
Garoufalidis who presented concrete problems requiring these techniques and
 R. D. Costin, S. Tanveer and J. McNeal for very useful discussions. Work supported in part
by NSF grants    DMS-0600369 and DMS-1108794.

\end{document}